\documentclass{amsart}
\usepackage{latexsym}

\usepackage{amsmath}
\usepackage{amsfonts}
\usepackage{amssymb}
\usepackage{eucal}
\title{How Dennis realized  he had `invented'  $L_\infty$- algebras a.k.a. strongly homotopy Lie algebras}
\author{jim stasheff}

\begin{document}

\maketitle
{\begin{abstract}

\large{This note is an attempt to rediscover
how Dennis realized that he had discovered/invented $L_\infty$- algebras.}
\end{abstract}

\baselineskip20pt

\section{Introduction}
\quad \quad $L_\infty$- algebras appeared independently  in the mid-1980's in a
supporting role in deformation theory in my work with Mike Schlessinger  \cite{Sjds,SS,jds-ms}
on rational homotopy theory and in correspondence between  Drinfel'd and Schechtman 
\cite{Drinfeld:1983,Drinfeld-to-Schechtman}
as well as in mathematical physics. They were in fact implicit in Dennis' models.  Here is my attempt to rediscover
how Dennis realized what he had discovered.

\quad \quad The philosophy that every deformation
problem is controlled by a differential graded Lie algebra, but not uniquely even up to isomorphism, leads to consideration of morphisms `up to homotopy'  and hence to $L_\infty$ algebras (originally called `strongly homotopy Lie algebras' or, by Drinfel'd, `Lie-Sugawara algebras'.

\section{Schechtman-Drinfel'd-Hinich-Sullivan}

\subsection{Schechtman writes:}

\quad The starting point of our studies with Volodya Hinich \cite{hinich-schechtman} on algebraic higher
homotopy was the following example. To define the multiplication in his 
cohomology theory which he invented to construct higher regulators, Beilinson used
some multiplication of complexes which was commutative up to a homotopy.
One noticed that this multiplication was a part of a richer structure; namely
 BeilinsonÕs complexes turned out to possess a structure of a module over a
certain algebraic contractible operad (its terms were chain complexes of cubes).

Apparently Volodya Drinfel'd was also interested in these subjects. Drinfel'd wrote that
he had some ideas on the subject and asked
if I am interested in details. 
Here is the
beginning of his letter dated 09/28/83 \cite{Drinfeld:1983}:
{\quote \quad \quad ``I am reading your papers with Hinich you have sent to me. Both of them 
are
very well written. In connection with your ÃlittleÓ paper (ÓOn homotopy limit
of homotopy algebrasÓ)  \cite{hinich-schechtman} I have some questions...''}

The known letter from Drinfel'd I received in September 1988 \cite{Drinfeld-to-Schechtman}, just before my
first trip to USA. Among other things, we can find the following there:

{\quote \quad \quad ``A Lie-Sugawara DG-algebra is, by definition, a Z-graded space g plus a
degree 1 differential on the cofree cocommutative coalgebra generated by g with
the grading shifted below (the square of the differential is 0).''}

\vskip2ex \quad 
This is what is called an $L_\infty$-algebra now; of course Drinfeld gives also a definition of a ÓDG-Sugawara (co)commutativeÓ (now $C_\infty$) algebra.
Concerning our paper with Hinich in the Gelfand seminar, we have discussed
these subjects at IHES in the summer Õ90. Ginzburg and Kapranov were also
there; they conceived their famous paper during this summer.
At some point, it was recognized that the description of an $L_\infty$-algebra as a coderivation differential on a  cofree connected graded symmetric coalgebra identified the $L_\infty$-algebra
implicit in Sullivan's models (not necessarily minimal).

\subsection{Drinfel'd comments further:}
\quad{\quote \quad At the time of my correspondence with Schectmann I felt that I was trying
to understand something known rather than inventing new things. Maybe my
feeling was correct. I wouldn't be surprised if everything from my
correspondence with Schechtmann is already in Quillen's article "Rational
Homotopy Theory" (maybe except the word "operad"). }
\subsection{Vladimir Hinich recalls:}

\quad{\quote\quad  I visited Schechtman at Stony Brook in 1992 during my last PhD year at Weizmann
Institute. If I remember correctly, he told me about Drinfeld's  letter to him (written in 1988),
and this was, I think, the only "source" of our paper in Gel'fand volume.}
 
\subsection{Sullivan recalls:}

\quad\begin{quote}\quad  I invited Schechtman to CUNY in the early nineties because i had heard he had come upon a notion of a global Lie algebra up to homotopy,
motivated by work on deforming algebraic varieties where a  sheaf up to homotopy of dgLie algebras controlled the deformation theory.

\quad My motivation was trying to discretize the pde for fluid motion while preserving all known properties, energy conservation, helicity conservation, vorticity frozen in the fluid, etc.
 I had done this except that my discrete version of volume preserving vector fields 
had a bracket which satisfied Jacobi only up to homotopy. I wondered if there was an analogue in  Lie algebra of what you had done for  associative H spaces....

\quad I was struck by lightning when Schechtman revealed that his Lie algebra up to homotopy was nothing but a differential-derivation on a free graded commutative algebra.

\quad Then it was clear that the various forms of rational homotopy theory could be viewed as infinity versions of structures here and there:
the dgc infinity coalgebra on chains computing homology was Quillen's differential on the free Lie algebra,
the dgLie infinity structure on a Moore complerx computing homotopy was Quillen's pre-dual coalgebra of the free dgc algebra models coming from forms.

\quad To summarize: what was new for me was the  familiar structures of rational homotopy theory were just infinity versions of appropriate structures on chain complexes...
\end{quote}


\begin{thebibliography}{10}

\bibitem{Drinfeld:1983}
V.G. Drinfeld, \emph{Letter to {S}chechtman}, 1983.

\bibitem{Drinfeld-to-Schechtman}
\bysame, \emph{Letter to {S}chechtman}, 1988.

\bibitem{hinich-schechtman}
V.~A. Khinich and V.~V. Shekhtman, \emph{The homotopy limit of homotopy algebras}, Uspekhi Mat. Nauk \textbf{41} (1986), no.~3(249), 205--206. \MR{854263}

\bibitem{jds-ms}
M.~Schlessinger and J.~Stasheff, \emph{Deformation theory and rational homotopy type}, {\tt{arXiv:1211.1647}}.

\bibitem{SS}
M.~Schlessinger and J.~D. Stasheff, \emph{The {L}ie algebra structure of tangent cohomology and de\-for\-ma\-tion theory}, J. of Pure and Appl. Alg. \textbf{38} (1985), 313--322.

\bibitem{Sjds}
J.~Stasheff, \emph{Rational homotopy theory -- obstructions and deformations}, Proc. Conf. on Algebraic Topology, Vancouver, 1977, LNM 673, pp.~7--31.

\end{thebibliography}

\end{document}